\documentclass[11pt]{amsart}
\makeatletter

\title{  Cohen-Macaulay Fiber Cones}
\author{ Clare D'Cruz}
\address
{Chennai Mathematical Institute, G. N. Chetty Road, T. Nagar, Chennai
600 017 India}
\email{clare@@cmi.ac.in}
\author{K. N. Raghavan} 
\address
{Chennai Mathematical Institute, G. N. Chetty Road, T. Nagar, Chennai
600 017 India}
\email{knr@@cmi.ac.in} 
\author{J. K. Verma}
\address
{Indian Institute of Technology Bombay Powai, Mumbai 400 076 India}
\email{jkv@@math.iitb.ac.in}
\newcommand{\ncom}{\newcommand}
\ncom{\beqn}{\begin{eqnarray*}}
\ncom{\eeqn}{\end{eqnarray*}}
\ncom{\beq}{\begin{eqnarray}}
\ncom{\eeq}{\end{eqnarray}}
\ncom{\been}{\begin{enumerate}}
\ncom{\eeen}{\end{enumerate}}
\ncom{\nno}{\nonumber}
\ncom{\hs}{\mbox{\hspace{.25cm}}}
\ncom{\rar}{\rightarrow}
\ncom{\lrar}{\longrightarrow}
\ncom{\Rar}{\Rightarrow}
\ncom{\noin}{\noindent}
\newtheorem{thm}{Theorem}[section]
\newtheorem{lemma}[thm]{Lemma}
\newtheorem{cor}[thm]{Corollary}
\newtheorem{pro}[thm]{Proposition}
\newtheorem{example}[thm]{Example}
\newtheorem{definition}[thm]{Definition}
\newtheorem{remark}[thm]{Remark}
\ncom{\bt}{\begin{thm}}
\ncom{\et}{\end{thm}}
\ncom{\bl}{\begin{lemma}}
\ncom{\el}{\end{lemma}}
\ncom{\bco}{\begin{cor}}
\ncom{\eco}{\end{cor}}
\ncom{\bp}{\begin{pro}}
\ncom{\ep}{\end{pro}}
\ncom{\bex}{\begin{example}}
\ncom{\eex}{\end{example}}
\ncom{\bd}{\begin{definition}}
\ncom{\ed}{\end{definition}}
\ncom{\brm}{\begin{remark}}
\ncom{\erm}{\end{remark}}
\ncom{\comx}{I\!\!\!\!C}
\ncom{\proj}{I\!\!\!P}
\ncom{\zee}{$Z\!\!\!\!Z$}
\ncom{\ze}{Z\!\!\!\!Z}
\ncom{\Q}{$I\!\!\!\!Q$}
\ncom{\N}{I\!\!N}
\ncom{\sz}{\scriptsize}
\ncom{\CM}{Cohen-Macaulay }
\ncom{\sop}{system of parameters}
\ncom{\eop}{\hfill{$\Box$}}
\ncom{\tfae}{the following are equivalent:}
\ncom{\mm}{minimal multiplicity }
\ncom{\f}{\frac}
\ncom{\la}{\lambda}
\ncom{\si}{\sigma}
\ncom{\ssize}{\scriptsize}
\ncom{\al}{\alpha}
\ncom{\be}{\beta}
\ncom{\Si}{\Sigma}
\ncom{\ga}{\gamma}
\ncom{\kbar}{\overline{\kappa}}
\ncom{\bib}{\bibitem}
\ncom{\sst}{\subset}
\ncom{\sms}{\setminus}
\ncom{\seq}{\subseteq}
\ncom{\est}{\emptyset}
\ncom{\pf}{{\bf Proof: }}
\ncom{\bighs}{\hspace{.5 cm}}
\ncom{\ulin}{\underline}
\ncom{\olin}{\overline}
\ncom{\bip}{\bigoplus}
\ncom{\sta}{\stackrel}
\begin{document}

\maketitle

\section{Introduction}

Let $(R,m)$ be a local ring and let $I$ be an ideal of $R$.
The {\em fiber cone} of $I$ is the graded ring $F(I):= \oplus_{n \geq
0} I^n / mI^n$. In this paper we will give a necessary and
sufficient condition, in terms of Hilbert series of $F(I)$,  for it
to be \CM (CM). We apply this criterion to a variety of
examples.

\CM fiber cones have been studied by K. Shah \cite{shah} and 
T.~~Cortadellas and S.~~Zarzuela \cite{cor}. Shah characterized \CM
fiber  cones  terms of  their
Hilbert function. To recall his result let us
introduce some notation and terminology. An ideal $J \seq I$ is
called a reduction of $I$ if $JI^{n} = I^{n+1}$ for some $n$. We
say $J$ is a minimal reduction of $I$ if whenever $K \seq I$ and
$K$ is a reduction of $I$, then $K=J$. These concepts were
introduced by Northcott and Rees in [NR]. They
showed that when $R/m$ is infinite, any minimal reduction of
$I$ is minimally generated by $\dim \,F(I)$ elements. The Krull
dimension of the fiber cone $F(I)$ is an important invariant of
$I$ and it is called the {\em analytic spread} of $I$. We denote it by
$a(I)$. It is known \cite{rees2} that ht$\; I \leq a(I) \leq
\min ( \mu(I), \dim \, R)$ where $\mu(I) = \dim \, I/mI$. Let
$\ell$ denote the length of the module.  We can
now state the necessary condition found by Shah for $F(I)$ to be
CM. 

\bt
\label{shah}
{\em (Shah)}
Let $(R,m)$ be a local ring and let $I$ be an ideal in $R$ of
analytic spread $a$. Let $J$ be a minimal
reduction of $I$. If $F(I)$ is CM, then for all $n \geq 0$
$$
\mu (I^n) 
= \sum_{i \geq 0} 
  \left[ \mu(I^i) - \ell \left(\f{JI^{i-1}}{JI^{i-1} \cap mI^i}
                                  \right) \right] 
      {n+a-i-1 \choose a-1} .
$$ 
 
\et

\noin

K.~Shah conjectured that the converse is  also true. 
In this paper we prove the converse.  The proof
is presented in section two. In section three we define ideals
of minimal mixed multiplicity in CM local rings and calculate
the Hilbert series of their fiber cone. It turns out that for
such ideals the fiber cone is CM if and only if $r(I) \leq 1$.
By using this result we construct zero dimensional ideals in a
two dimensional regular local ring whose fiber cones are not CM.
Section four is about fiber cones of ideals generated by
quadratic sequences.
In the final section we provide a few other examples of CM fiber
cones.  

{\bf Acknowledgements:} The first author wishes to thank the National Board of Higher
Mathematics for financial support and the Indian Institute of
Technology, Bombay, where the research on this paper was carried
out. 
The second named author wishes to acknowledge financial help 
in the form a grant from the
Department of Science and Technology of the Government of India
to the SPIC Mathematical Institute  for purchase of
computing equipment on which calculations using {\it Macaulay}
could be
performed.    The last named author wishes to acknowledge
the hospitality of the SPIC Mathematical Institute during a visit
on which part of the work on this paper was done.

\section{\CM fiber cones}

In this section we
state K.~Shah's result in terms of Hilbert
series of $F(I)$. For the sake of completeness we give an
alternate proof and prove the converse. 
We will denote the multiplicity of $F(I)$ by $e(F(I))$. We know
that if $F(I)$ is CM then 
$e(F(I)) \geq \mu(I) - a + 1$ \cite{abh}. When
equality holds, we say that $F(I)$ has { \em minimal
multiplicity}.  The reduction number of $I$ with respect to a
minimal reduction $J$ of $I$ is defined as 
$$ 
r_J (I) = \min\{ n \geq 0 \;|\; JI^n =I^{n+1} \}.  
$$ 
The reduction number of $I$ denoted by $r(I)$ is defined to be
the minimum of $r_J(I)$ where $J$ varies over all minimal
reductions of $I$. 

\bt
\label{main}
Let $(R,m)$ be a local ring and let $I$ be any ideal of $R$.
Let $J$ be any minimal reduction of $I$. Then the following are
equivalent:
\been
\item
 $F(I)$ is CM.
\item
 $ \displaystyle{
H(F(I),t) 
= \f{1}{(1 -t)^{a}}  \sum_{i=0}^{r} \ell \left(
                   \f{I^i}{JI^{i-1}+mI^i} \right) \, t^i}
$
where $r=r_J(I)$ and $a= \dim \, F(I)$.
\item
 $ \displaystyle{ e(F(I)) = \sum_{n=0}^{r} \ell \left(
\f{I^{n}}{JI^{n-1}+mI^n} \right)}$.
\eeen
\et
\pf $(1) \Rar (2)$. Suppose that $F(I)$ is CM. Since $JF(I)$ is
generated by a homogeneous system of parameters, it is generated
by a regular sequence. Thus 
$$ 
H(F(I),t) 
= \f{1}{(1 -t)^{a}} H\left(\f{F(I)}{JF(I)},t \right) .
$$
Notice that 
$$
JF(I)
= \bip_{i=1}^{r} \left( \f{JI^{i-1} + mI^i}{mI^i} \right)
   \bip \left[\bip_{i=r+1}^{\infty} 
   \left( \f{I^i}{mI^i}\right) \right].
$$
Hence 
$$
\f{F(I)}{JF(I)}
= \f{R}{m}  \bip \left[
  \bip_{i=1}^{r} \left( \f{I^i}{JI^{i-1} + mI^i} \right) \right].
$$
This gives
$$
H\left(\f{F(I)}{JF(I)},t \right) 
= \sum_{i=0}^{r} \ell \left( \f{I^i}{JI^{i-1} + mI^i} \right)  t^i.
$$
which gives $(2)$. 

$(2) \Rar (3)$. Put $t=1$ in the numerator of $H(F(I),t)$.

$(3) \Rar (1)$. Let ${\mathcal M}$ denote the maximal homogeneous
ideal of $F(I)$. It is enough to show that $F(I)_{\mathcal M}$ is CM
by \cite{mr}. Since $F(I)$ is homogeneous graded ring, the
associated graded ring $gr({\mathcal M }F(I)_{\mathcal M}, F(I)_{\mathcal 
M}) \cong
F(I)$. Thus $e(F(I)) = e(F(I)_{\mathcal M})$.
It is also clear that 
$\ell ( F(I)/JF(I)) = \ell (F(I)_{\mathcal M} / JF(I)_{\mathcal M})$.
Thus $F(I)_{\mathcal M}$ is CM and hence $F(I)$ is so.

\brm\label{remark}
{\em Let $(R,m)$ be a local ring and let $I$ be any ideal of $R$.
If $F(I)$ is \CM, then by theorem~\ref{main}, for any minimal
reduction $J$ of $I$, $r_{J}(I)$ is the 
degree of the numerator of the Hilbert series of $F(I)$ in the
reduced form. Hence the reduction number $r_{J}(I)$ is independent
of the choice of the minimal reduction $J$ of $I$. This
observation is also proved, by different means, by Huckaba and
Marley  in \cite{hm}}.
\erm

\bco
\label{main1}
{\em (Shah)}
Let $(R,m)$ be a local ring. Suppose $I$ is an ideal of $R$
of analytic spread $a$ having a minimal reduction $J= (x_1,
\ldots, x_a)$ with $JI=I^2$. If $x_{1}, \ldots, x_{a}$ is a
regular sequence, then $F(I)$ is CM with minimal multiplicity. 
\eco
\pf We can extend $x_1, \ldots, x_a$ to a minimal basis of $I$.
So let $\mu(I) = a+p$ and $I= (x_1, \ldots, x_a, y_1, \ldots,
y_p)$. Set $K=(y_1, \ldots, y_p)$. Then 
$I^n = (J^n, J^{n-1}K, \ldots, K^n)$. For $i \geq 2$, 
$J^{n-i} K^i \seq J^{n-i} I^i = I^n = J^{n-1}I$ since $JI =
I^2$. Thus $I^n = (J^n, J^{n-1}K)$. Let $M_1, \ldots, M_r$ be
all the monomials of degree $n$ in $x_1, \ldots, x_a$ and let $m_1,
\ldots, m_N$ be 
all the monomials in $x_1, \ldots, x_a$ of degree $n-1$. We
claim that $I^n$ is minimally generated by the elements
$$
\{M_1, \ldots, M_r\}
\cup \{m_i y_j \, | \, i=1, \ldots, N; \, j=1, \ldots, p\}.
$$
Indeed, let them satisfy the relation 
$$
  \sum_{i=1}^{r} a_{i} M_i 
+ \sum_{j=1}^{p} \sum_{i=1}^{N}  b_{ij} m_i y_j  = 0.
$$
Let $b_{\alpha \beta} \not \in m$. Rewrite the above relation as
\beqn
  \sum_{i=1}^{r} a_{i} M_i 
+ \sum_{j=1}^{p} 
  \sum_{\sta{i=1}{i \not = \alpha}}^{N} b_{ij} m_i y_j 
+  m_{\alpha} \sum_{j=1}^{p} b_{\alpha j} y_{j}= 0.
\eeqn
Hence 
\beqn
\sum_{j=1}^{p} b_{\alpha j} y_{j}
\in (J^{n}, m_1, \ldots, \hat{m}_{\alpha}, \ldots, m_{N}) : m_{\alpha}
\eeqn
which gives
$$
y_{\beta} \in (J, y_1, \ldots, \hat{y}_{\beta}, \ldots, y_{p})
$$
which gives a contradiction. Thus 
\beqn
H(F(I), t) 
&=& \sum_{n=0}^{\infty}\mu(I^n) t^n \\
&=& \sum_{n=0}^{\infty} 
    \left[ {a-1+n \choose a-1} + {a-2+n \choose a-1} p \right] t^n \\ 
&=& \sum_{n=0}^{\infty}  {a-1+n \choose a-1} t^n
   + p \sum_{n=0}^{\infty} 
     \left[ {a-1+n \choose a-1} 
  - {a-2+n \choose a-2} \right] t^n \\ 
&=& \f{1}{(1-t)^{a}} 
   + p \left[ \f{1}{(1-t)^{a}} - \f{1}{(1-t)^{a-1}} \right]\\
&=& \f{1+pt}{(1-t)^{a}}
\eeqn
It remains to show that $p = \ell (I/J+mI)$. In fact 
$$
\ell \left(\f{I}{J+mI} \right)
= \ell \left(\f{I}{mI} \right) - \ell \left(\f{J+mI}{mI} \right)
= \mu(I) - \ell \left(\f{J}{mJ} \right)
= \mu(I) - \mu(J)
= p..
$$
Therefore $F(I)$ is CM. Since $e(F(I)) = 1 + \mu(I) - a$, $F(I)$
is CM with minimal multiplicity. 

\section{Fiber cones of ideals with minimal mixed multiplicity}

Let $(R,m)$ be a local ring. Define the function $H: \N^{g} \rar \N$
for $g$  $m$-primary ideals $I_1, \ldots, I_g$ as follows:
\beqn
H(r_1, \ldots, r_g)  
= \ell \left( \f{R}{I_1^{r_1} \cdots I_g^{r_g}} \right).
\eeqn
It was proved in \cite{teis} that $H(r_1, \ldots, r_g)$ is given
by a polynomial $P(r_1, \ldots, r_g)$ with rational coefficients
for all large $r_1, \ldots, r_g$. The total degree of $P$ is $d$
and the terms of degree $d$ in $P$ can be written as
\beqn
\sum_{i_1+ \cdots+i_g\,=\,d} e(I_1^{[i_1]} | \cdots | I_g^{[i_g]})
{r_1 + i_1 \choose i_1} \cdots {r_g + i_g \choose i_g}
\eeqn
for certain positive integers $e(I_1^{[i_1]} | \cdots | I_g^{[i_g]})$.
These positive integers are called mixed multiplicities of $I_1,
\ldots, I_g$. When $g=2$ we adopt the notation 
\beqn
e_i(I_1|I_2) = e(I_1^{[d-i]} | I_2^{[i]}).
\eeqn
Rees \cite{rees1} showed that $e_0(I_1|I_2)=e(I_1)$ and
$e_d(I_1|I_2)=e(I_2)$. Other mixed multiplicities can also be
interpreted as multiplicities of certain parameter ideals.
Teissier and Risler \cite{teis} provided an interpretation of
mixed multiplicities in terms of superficial elements. Rees
\cite{rees4} introduced joint reductions to study them. 
A set of elements $x_1, \ldots, x_d$ is a called a {\em joint
reduction} of $I_1, \ldots, I_d$ if $x_i \in I_i$ for $i=1,
\ldots, d$ and there exists a positive integer $n$ so that 
\beqn
  \left[ \sum_{j=1}^{d} x_j \, I_1 \cdots \hat{I_j} \cdots I_d \right] 
   (I_1 \cdots I_d)^{n-1}
=  (I_1 \cdots I_d)^n.
\eeqn
Rees proved that if $R/m$ is infinite, then joint reductions
exist and 
\beqn
e(I_1^{[i_1]} | \cdots | I_g^{[i_g]}) = e(x_1, \ldots, x_d)
\eeqn
where $(x_1, \ldots, x_d)$ is any joint reduction of the set of
ideals 
$I_1, \ldots, I_1,  \ldots, I_g, \ldots, I_g$
where $I_j$ is repeated $i_j$ times for $j=1, \ldots, g$.. 

Abhyankar \cite{abh} showed that for a CM local ring $(R,m)$, $\mu(m) -
\dim \, R + 1 \leq e(R)$. We extend this to $m$-primary ideals. 

\bl
Let $(R,m)$ be a CM local ring of dimension $d$. Then 
$$
e_{d-1}(m|I) \geq \mu(I) - d + 1.
$$
\el
\pf By passing to the ring $R(X)=R[X]_{mR[X]}$ we may assume
that $R/m$ is infinite. Let $(x, a_1, \ldots, a_{d-1})$ be a
joint reduction of $m, I, \ldots, I$ where $I$ is repeated $d-1$
times. Consider the map 
\beqn
f: \f{R}{I} \bip \left( \f{R}{m} \right)^{d-1} 
\, \lrar \, \f{(x, a_1, \ldots, a_{d-1})}
              {xI+(a_1, \ldots, a_{d-1})m}
\eeqn
defined as $f(a', b_1', \ldots, b_{d-1}') =
(xa + a_1b_1 + \cdots + a_{d-1}b_{d-1})^{\,\prime}$ 
where primes denote the residue classes. Let
$K=\mbox{Kernel}\,f$. Then 
\beqn
\ell(K) + \ell \left( \f{R}{xI+(a_1, \ldots, a_{d-1})m} \right)
  - e_{d-1}(m|I) 
= \ell \left( \f{R}{I} \right) + d-1.
\eeqn
Hence
\beqn
    e_{d-1}(m|I) 
&=& \ell(K) + \ell \left( \f{R}{Im} \right) 
  - \ell \left( \f{R}{I} \right)
  - d + 1  - \ell \left( \f{R}{Im} \right) 
  + \ell \left( \f{R}{xI +(a_1, \ldots, a_{d-1})m} \right)\\
&=& \ell(K) + \mu(I) - d + 1 
  + \ell \left( \f{Im}{xI +(a_1, \ldots, a_{d-1})m} \right).
\eeqn
Hence
$$
e_{d-1}(m|I) \geq \mu(I) - d + 1.
$$

\bco
\label{}
{\em [cf. Sa, pg 49]}
Let $(R,m)$ be a one-dimensional CM local ring. Then for any
$m$-primary ideal $I$ of $R$
$$
\mu(I) \leq e(R).
$$
\eco
\pf When $d=1$, $e_{d-1}(m|I) = e_0(m|I) = e(R)$. 

\bco
Let $(R,m)$ be a $d$-dimensional CM local ring. Then 
$$
e(m) \geq \mu(m) - d + 1.
$$
\eco
\pf Put $I=m$ and use the fact that $e_i(m|m)=e(m)$ for all $i$.

\bd
{\em Let $(R,m)$ be a CM local ring, $I$ an $m$-primary ideal of $R$. We
say that $I$ has {\em minimal mixed multiplicity} if }
$$
e_{d-1}(m|I) = \mu(I) - d + 1.
$$
\ed

We now find the Hilbert series of $F(I)$ when $I$ is an
$m$-primary ideal with minimal mixed multiplicity in a
$d$-dimensional CM local ring. We begin with one dimensional
local rings.  

\bp
Let $(R,m)$ be a CM local ring of dimension one with $\mu(I) =
e(R) :=e$ for an $m$-primary ideal $I$ of $R$. Then 
\beqn
H(F(I),t) = \f{1 + (e-1)t}{1-t}.
\eeqn
\ep
\pf By passing to the ring $R(X)$ we may assume that $R/m$ is
infinite. Let $xR$ be a minimal reduction of $m$. Then for any
$n \geq 1$, 
\beqn
e(R) = \ell \left( \f{R}{xR} \right) 
     = \ell \left( \f{R}{xR} \right) 
     + \ell \left( \f{xR}{xI^n} \right) 
     - \ell \left( \f{R}{I^n} \right)
     = \ell \left( \f{I^n}{xI^n} \right) .
\eeqn
Since $\ell(I/mI) = e(R) = \ell (I/xI)$, it follows that $mI =
xI$. Thus $xI^n = mI^n$ which gives that $\mu(I^n) = e(R)$ 
for all $n \geq 1$. Hence 
\beqn
H(F(I), t) = 1 + \sum_{n=1}^{\infty} e \,t^n
           = 1 - e + \f{e}{1-t}
           = \f{1 + (e-1)t}{1-t}.
\eeqn

\bt
\label{depth}
Let $(R,m)$ CM local ring of dimension $d \geq 2$.  Let $I$ be
an $m$-primary ideal having minimal mixed multiplicity.  Then  
\beqn
H(F(I), t) = \f{1 + (\mu(I) - d)t}{(1-t)^d}
\eeqn
\et
\pf By passing to the ring $R(X)$ we may assume that $R/m$ is
infinite. Let $x, a_1, \ldots, a_{d-1}$ be a joint reduction of
$(m,I, \ldots, I)$ where $I$ is repeated $d-1$ times. Then
$e_{d-1}(m|I) = \ell (R/(x,a_{1}, \ldots, a_{d-1}))$ and $\mu(I) =
d-1 + \ell (R/(x,a_{1}, \ldots, a_{d-1}))$. We 
will show that for all $n \geq 1$,
\beqn
\mu(I^n) = {n + d-2 \choose d-2} + {n +  d-2 \choose d-1} e_{d-1}(m|I). 
\eeqn
Fix $n$. Put $r = {n + d-2 \choose d-2}$. Let $M_{1}, \ldots,
M_{r}$ denote all the monomials of degree $n$ in 
$a_{1}, \ldots, a_{d-1}$. 
Consider the map 
\beqn
f: \f{R}{I^n} \oplus \left( \f{R}{m} \right)^{r}
\, \lrar \, \f{(x, (a_{1}, \ldots, a_{d-1})^n)}
              {xI^n + (a_{1}, \ldots, a_{d-1})^n m}  
\eeqn
given by 
$f(y',z_{1}', \ldots, z_{r}') 
= (yx + \sum_{i=1}^{r} z_{i} M_{i})'$ 
where primes denote the residue classes. If 
$f(y',z_{1}', \ldots, z_{r}') = 0$, then
$yx + \sum_{i=1}^{r} z_{i} M_{i}
\in xI^n + (a_{1}, \ldots, a_{d-1})^n m$. Hence there exist $b
\in I^n$ and $c_{1}, \ldots, c_{r} \in m$ such 
that 
\beq
\label{mmm2}
yx + \sum_{i=1}^{r} z_{i} M_{i}  
= xb + \sum_{i=1}^{r} c_{i} M_{i} .
\eeq
Thus $x(y-b) = \sum_{i=1}^{r} M_{i} (c_{i} - z_{i})$. Since
$(x, a_{1}, \ldots, a_{d-1})$ is a regular 
sequence, 
$y - b \in (M_{1}, \ldots, M_{r})
       = (a_{1}, \ldots, a_{d-1})^n$. 
Hence 
$y = b + \sum_{i=1}^{r} d_{i} M_{i}$ 
for some $d_{1}, \ldots, d_{r} \in R$. Substituting for $y$ in
(\ref{mmm2}), we get  
$\sum_{i=1}^{r} M_{i} (z_{i}  - c_{i}  +  xd_{i}) = 0$. 
Since $a_{1}, \ldots, a_{d-1}$ is a regular sequence, 
$z_{i} - c_{i} + xd_{i} \in (a_{1}, \ldots, a_{d-1})$ for all
$i=1, \ldots, r$. Hence $z_{i} \in m$ for all $i=1, \ldots, r$.
Thus $f$ is an isomorphism which implies that 
\beqn
   \ell \left( \f{R}{xI^n + (a_{1}, \ldots, a_{d-1})^n m} \right) 
-  \ell \left( \f{R}{(x, (a_{1}, \ldots, a_{d-1})^n)} \right) 
=  \ell \left( \f{R}{I^n} \right)  + r.
\eeqn
Put $n = 1$ to get 
\beqn
   \ell \left( \f{R}{xI + (a_{1}, \ldots, a_{d-1}) m} \right)  
 - \ell \left( \f{R}{I} \right) 
=   r + e_{d-1}(m|I) 
=  \mu(I) 
=  \ell \left( \f{R}{mI} \right) 
 - \ell \left( \f{R}{I} \right).
\eeqn
Hence $Im = xI + (a_{1}, \ldots, a_{d-1}) m$. Thus for all $n \geq
1$, $I^n m = xI^n + (a_{1},\ldots, a_{d-1})^n m$ by an easy
induction on $n$. Hence for all $n \geq 1$, 
\beqn
\mu(I^n) = r + {n + d-2 \choose d-1}  e_{d-1}(m|I).
\eeqn
Therefore
\beqn
H(F(I), t) 
&=& \sum_{n=0}^{\infty} 
    \left[r + {n+d-2 \choose d-1} e_{d-1}(m|I) \right]\,t^n \\ 
&=& \sum_{n=0}^{\infty} 
    \left[r (1 - e_{d-1}(m|I)) 
    + {n+d-1 \choose d-1} e_{d-1}(m|I) \right]\,t^n \\ 
&=& \f{1 - e_{d-1}(m|I)}{(1-t)^{d-1}} 
  + \f{e_{d-1}(m|I)}{(1-t)^d} \\
&=& \f{1 + [e_{d-1}(m|I) - 1]t}{(1-t)^d} \\
&=& \f{1 + [\mu(I) - d]t}{(1-t)^d}.
\eeqn

\bt
\label{cor}
Let $I$ be an $m$-primary ideal in a CM local ring $(R,m)$.
Suppose $I$ has minimal mixed multiplicity. Then $F(I)$ is CM
if and only if $r(I) \leq 1$. 
\et
\pf If $r(I) \leq 1$, then $F(I)$ is CM by
corollary~\ref{main1}. If $F(I)$ is CM, then 
$$
r(I) = \mbox{ the degree of the numerator of } H(F(I),t) 
       \leq 1. 
$$

We now consider fiber rings of contracted ideals in  two
dimensional regular local rings. An ideal $I$ of a two
dimensional regular local ring $(R,m)$ is called {\em contracted} if
there is an $x \in m \sms m^2$ such that $IR[m/x] \cap R =I$
[ZS, App.~5]. The order of $I$, $o(I)$ is by definition
max$\{n \, | \, I \seq m^n \}$. I is  contracted if and
only if $\mu(I) = 1 + o(I)$ [HS, Theorem~2.1]. If $I$ is
contracted from $S=R[m/x]$, then for all $n$, $I^n$ is also
contracted from $S$. Any $m$-primary integrally closed ideal is
contracted [ZS, App.~5]. 

\bco
\label{cont}
Let $(R,m)$ be a two dimensional regular local ring. Let $I$ be
an $m$-primary contracted ideal of order $r$. Then 
$$
H(F(I), t) = \f{1 + (r-1)t}{(1-t)^2}.
$$
\eco
\pf By [V1, Theorem~4.1], $e_1(m|I)=r = o(I)$. Since $I$ is an ideal
of minimal mixed multiplicity, we can apply
Theorem~\ref{depth}.

In view of the above theorem~\ref{cor}, to produce a non CM
fiber cone it is enough to find an ideal having minimal mixed
multiplicity whose reduction number exceeds one. We do so in the
next example. 

\bex
\label{noncm}
{\em
We construct  contracted ideals of reduction number greater than
one which have non CM fiber cones. Consider the ideals
\beqn
I = (x^r, xy^{r+1}(x,y)^{r-2}, y^{2r+1}).
\eeqn
of $k[[x,y]]$ where $k$ is a field. 
Since $\mu(I) = r+1 = 1 + o(I)$, $I$ is a contracted ideal and
by corollary~\ref{cont}, 
$$
H(F(I), t) = \f{1 + (r-1)t}{(1-t)^2}.
$$
If $F(I)$ were CM, then by corollary~\ref{cor}, $r(I) \leq 1$.
We now show that $r(I) >1$ for $r \geq 3$. By [V2,
Theorem~3.2], $r(I) \leq 1$ if and only if $e(I) = \ell(R/I^2) - 2 \ell
(R/I)$. Since $I$ is integral over $(x^r, y^{2r+1})$, $(x^r,
y^{2r+1})$ is a minimal reduction of $I$ and $e(I) = r(2r+1)$.
Now 
\beqn
I^2
= ( x^{2r}, x^{r+1}y^{r+1}(x,y)^{r-2}, x^r y^{2r+1}, 
    x^2 y^{3r+1}(x,y)^{r-3}, x y^{4r}, y^{4r+2}).
\eeqn
Therefore 
\beqn
   \ell \left( \f{R}{I^2} \right)
&=&  \sum_{i=0}^{r-2} (r + 1 + i) + 2r + 1 
    + \sum_{i=0}^{r-3}(3r+1+i) + 4r + 4r+2\\
&=&   (r+ 1)(r-1) + {r-1 \choose 2} + 2r + 1
    + (3r + 1)(r-2) + {r-2 \choose 2} + 8r + 2\\
&=&   {r-1 \choose 2} + {r-2 \choose 2} + 4r^{2} + 5r
\eeqn
and
\beqn
\ell \left( \f{R}{I} \right)
&=& \sum_{i=0}^{r-2} (r + 1 + i) + 2r + 1 \\
&=& (r+ 1)(r-1) + {r-1 \choose 2} + 2r + 1\\
&=&  {r-1 \choose 2} + r^2 + 2r.
\eeqn
Therefore if $r \geq 3$, then
\beqn
\ell \left( \f{R}{I^2} \right) - 2 \ell \left( \f{R}{I} \right)
&=& 4r^2 + 5r  + {r-2 \choose 2} 
-  2\left[r^2 + 2r + {r-1 \choose 2} \right] \\
&=& 2r^2 + r - (r-2)\\
&<& 2r^2  + 2r\\
&=& e(I)
\eeqn}
\eex

\bex {\em
\label{}
[cf. [HJLS, (6.3)],[HM, Corollary (3.5) and Example (3.9)]
There exist examples of $m$-primary ideals in $R=k[[x,y]]$ of
arbitrary reduction number and whose fiber ring is CM. Let
$I=(x^r, x^{r-1}y, y^r)$, $r \geq 1$. Then $F(I)$ is CM for
all $r \geq 1$.  
Since $I$ is generated by homogeneous polynomials of equal
degree in $R$, $F(I) \cong k[x^r, x^{r-1}y, y^r]$. 
It is easy to see that $k[x^r, x^{r-1}y, y^r] \cong
k[t_1,t_2,t_3]/(t_2^{r} -t_1^{r-1}t_3)$. Hence $F(I)$ is CM and
its Hilbert series is 
$$
H(F(I),t) = \f{1+t+ \cdots + t^{r-1}}{(1-t)^2}.
$$
It follows from theorem~\ref{main} that the reduction number of
$I$ is independent of the reduction chosen and is equal to
$r-1$.}
\eex

\bex
{\em 
Let $R = k[[x,y,z]]$, $I= (x^{3}, y^{3}, z^{3}, xy, xz, yz)$.
This ideal was studied by Huneke and Lipman \cite{lipman}. We
claim that $I$ has minimal mixed multiplicity. 
It is easy to see that
\been
\item
$J= (x^{3}+yz, y^{3}+z^{3}+xz, xz+xy)$
is a  minimal reduction of $I$ and $JI=I^{2}$. 

\item
$(yz, y+z, x)$ is a joint reduction of the set of
ideals (I, m, m) and 
$Im^{2} = yz m^{2} + (y+z) Im + x Im$. 

\item
$(yz, xy+xz, x+y+z)$ is a joint reduction of the set of
ideals $(I, I, m)$ and 
$I^{2}m  = yz Im + (xy+xz)Im + (x+y+z) I^{2}$.
\eeen
Hence 
\beqn
e(I)       &=& e(J) = 11\\
e_{1}(m|I) &=& e(yz, y+z, x)=2 \\
e_{2}(m|I) &=& e(yz, xy+xz, x+y+z) = 4.
\eeqn
Since $\mu(I) = 6 = e_{2}(m|I) + 2$, $I$ has minimal mixed
multiplicity. Hence by corollary~\ref{main1}, $F(I)$  is CM.
}
\eex

\input{amssym.def}

\newsymbol\lneq 230C
%\ncom{\lneq 230C }{\not \leq}
\renewcommand{\thefootnote}{}

\newcommand{\rk}{\rm rk}
\newcommand{\ul}{\underline}
\newcommand{\ve}{\varepsilon}
\newcommand{\supp}{\rm supp}
\newcommand{\bq}{{\Bbb Q}}
\newcommand{\bn}{{\Bbb N}}
\newcommand{\bz}{{\Bbb Z}}
\newcommand{\bc}{{\Bbb C}}
\newcommand{\cm}{{\mathcal M}}
\newcommand{\cn}{{\mathcal N}}
\newcommand{\cp}{{\mathcal P}}
\newcommand{\fg}{{\mathfrak g}}
\newcommand{\fs}{{\mathfrak s}}
\newcommand{\cf}{{\mathcal F}}

\newcommand{\scriptA}{{\mathcal A}}
\newcommand{\scriptF}{{\mathcal F}}
\newcommand{\gr}{{\rm gr}}
\newcommand{\osum}{\oplus}

\newcommand{\lr}{\longrightarrow}
\newcommand{\wt}{\widetilde}
\newcommand{\ol}{\overline}
\newcommand{\nk}{{\mathfrak A}}
\newcommand{\ds}{\displaystyle}
\newcommand{\bh}{{\Bbb H}}
\newcommand{\ci}{{\mathcal I}}
\newcommand{\tra}{\twoheadrightarrow}
\renewcommand{\char}{{\rm char}\,}

\section{Cohen-Macaulayness of fiber cones of ideals generated by
quadratic sequences.}

Let $R$ be a standard graded ring over a field,  let $M$ denote
the unique graded maximal ideal of $R$,  and let $I$ be an ideal 
generated by a homogeneous quadratic sequence.    Under some technical
assumptions on the quadratic sequence,
we show that a certain deformation of the fiber cone $F(I):=R[It]/MR[It]$
has a nice form.    
   From this we can deduce
the Cohen-Macaulayness of $F(I)$ and also a formula for its Hilbert
series.    (The Hilbert series can also be calculated 
from the formula in [RV] for the bigraded Hilbert series of
$R[It]$---just set $X=1$ in that formula.)
These results apply to the examples treated in [RV],  namely,
straightening-closed ideals in graded algebras with straightening law,
Huckaba-Huneke ideals of analytic deviation $1$ and $2$,  and Moral\`es-Simis
ideals defining the homogeneous co-ordinate rings of certain projective
space curves.

The main result of this section is stated in terms of stable linearizations
of quadratic sequences.
We recall the relevant definitions from [RS].
A subset $\Lambda$ of a finite
poset $(\Omega, \leq)$ is
an {\it ideal} if 
$$
\lambda \in \Lambda, ~~\omega \in \Omega, ~~\mbox{and}~~
\omega \leq \lambda ~~~~~~\Longrightarrow ~~~~~~~\omega \in \Lambda .
$$
If $\Lambda$ is an ideal of $\Omega$ and $\omega \in \Omega
\backslash \Lambda$ is such that $\lambda \in \Lambda$ for
every $\lambda\lneq \omega$, then $(\Lambda ,\omega
)$ is a {\it pair} of $\Omega$.  Given a set $\{ x_{\omega}
\mid \omega \in \Omega \}$ of elements of a ring $R$ and
$\Lambda \subseteq \Omega$, denote by $X_{\Lambda}$ the
ideal $(x_{\lambda} \mid \lambda \in \Lambda )$ of $R$
($X_{\Lambda} =0$ if $\Lambda$ is empty) and by $I$ the
ideal $X_{\Omega} = (x_{\omega} \mid \omega \in \Omega$).

\bd
\label{dquad} 
{\em
A set $\{ x_{\omega} \mid \omega \in \Omega \} \subseteq R$
is a {\it quadratic sequence} if for every pair $(\Lambda
,\omega )$ of $\Omega$ there exists an ideal $\Theta$ of
$\Omega$ such that
\begin{enumerate}
\item $(X_{\Lambda}: x_{\omega}) \cap I = X_{\Theta}$.
\item $x_{\omega} X_{\Theta} \subseteq X_{\Lambda} I$.
\end{enumerate}
Such an ideal $\Theta$ is said to be {\it associated} to the
pair $(\Lambda ,\omega )$.  This association need not be
unique---the set $\{ x_{\omega} \mid \omega \in \Omega \}$
of generators of $I$ may not be unshortenable---but
$X_{\Theta}$ is unique by 1.
}
\ed

\bd
\label{dlin}
{\em
A {\it linearization} of a poset
$\Omega$ of cardinality $n$ is a bijective map $\# : \Omega
\lr [1,n] :=\{1,\ldots ,n\}$ such that $\omega \leq \omega'
\Longrightarrow \# (\omega ) \leq \# (\omega ')$.
}
\ed

\noindent
Let $\{ x_{\omega} \mid \omega \in \Omega\}$ be a quadratic
sequence and $\# : \Omega \lr [1,n]$ be a fixed
linearization.  Identify $\Omega$ with $[1,n]$ via $\#$.
Then $([1,j-1],j)$ is a pair of $\Omega$ for every $j \in
[1,n]$.  Let
$$\begin{array}{rcl}
\Theta_j & = & \mbox{an ideal of $\Omega$ associated to
$([1,j-1],j)$} \\ [2mm]
I_j & = & ((x_1,\ldots ,x_{j-1}):x_j) \\ [2mm]
\Delta & = & \{ (j,k) \mid 1 \leq j \leq k \leq n,
~~x_j\cdot x_k \in (x_1,\ldots ,x_{j-1})\} .\\
\end{array}$$
For $k \in [1,n]$, set
$$\Psi_k = [1,k-1] \bigcup_{j \leq k \atop (j,k) \not\in
\Delta} \Theta_j ~~\mbox{and}~~\nk_k = X_{\Psi_k}.$$
Note that $\Psi_k$ is an ideal of $\Omega$ and that $\nk_k$
is independent of the choices of $\Theta_j$.

\bd
\label{dstab} 
{\em
A linearization $\# : \Omega
\lr [1,n]$ of the indexing poset $\Omega$ of a quadratic
sequence is {\it stable} if $I_k = (\nk_k:x_k)$ for every
$k$, $1\leq k\leq n$.  
}
\ed

We can now state the main result of this section:
\bp
Let $R$ be a standard graded algebra over a field $k$,  that is,
$R=\osum_{j\geq0}R_j=R_0[R_1]$ with $R_0=k$.    Let $\{x_\omega
| \omega\in\Omega\}\subseteq R$ be a quadratic sequence consisting
of homogeneous elements of $R$,  let $\#:\Omega\to\{1,\ldots,n\}$
be a stable linearization, and suppose that 
$$
\deg(x_1)\leq\cdots\leq\deg(x_n).
$$
Further assume that $x_1,\ldots,x_n$ form an unshortenable set of
generators of the ideal $I:=(x_1,\ldots,x_n)$,  that is,  $I\neq
(x_1,\ldots,\widehat{x_j},\ldots,x_n)$.    Let $M:=\osum_{j>0}R_j$
denote the irrelevant maximal ideal of $R$.   Then
$$
k[\Delta]:=k[T_1,\ldots,T_n]/(T_jT_k | (j,k)\in\Delta)
$$
is a deformation of the fiber cone $F(I):=R[It]/MR[It]$ of $I$.
In particular (1) $F(I)$ is Cohen-Macaulay if $k[\Delta]$ is so,
and (2) $F(I)$ has the same Hilbert series as $k[\Delta]$.
\ep

\noindent{\bf Proof.}
We borrow the notation and set-up of \S1 of [RS]
(see page~541).   Let $J$ denote the kernel of the
presentation map from the polynomial ring $\scriptA:=
R[T_1,\ldots,T_n]$ onto the Rees ring $R[It]$.   This map
is defined by $T_i\mapsto x_it$.    The presentation ideal of
the fiber cone $F(I)$ as a quotient of $\scriptA$ is then
$(J,M)$.    Let $\scriptF$ be the filtration on $\scriptA$ specified
on page~541 of [RS].     This filtration $\scriptF$ is
a refinement of the $(M,T_1,\ldots,T_n)$-adic filtration on 
$\scriptA$.    Theorem~1.4 of [RS] asserts that
$$
(*)~~~~~\gr_{\scriptF}(J)=(I_1T_1,\ldots,I_nT_n,T_jT_k| (j,k)\in\Delta),
$$
where $I_j:=((x_1,\ldots,x_{j-1}):x_j)$.    So it is natural to
expect the following:
$$
(\dag)~~~~~\gr_{\scriptF}(J,M)=(M,T_jT_k|(j,k)\in\Delta).
$$
The theorem clearly follows from ($\dag$) and we prove ($\dag$) below.

Since $M$ is homogeneous with respect to $\scriptF$ and $T_jT_k\in\gr_{
\scriptF}(J)$ from (*), it follows that $\gr_{\scriptF}(J,M)\supseteq
(M,T_jT_k|(j,k)\in\Delta)$.    To prove the other inclusion, it suffices
to prove the statement ($\ddag$) below.    Fix notation as in 3.2 of
[RS].      Let $F$ be a form of degree $d$ in $J$.  Let $Q$
be the subset of $P(n,d)$ such that $F=\sum_{p\in Q}a_pT^p$ where
$0\neq a_p\in R$ for $p\in Q$.   Assume further that the $a_p$ are 
homogeneous.   Let $q$ be the such element of $Q$ that
the initial form $F_*$ of $F$ with respect to $\scriptF$ equals
$a_qT^q$,
and let $s$ be the such element of $Q$ that $a_s\not\in M$ and 
the term $F_s=a_sT^s$
has the least degree in ${\N}_0^{n+1}$ among those with $a_p\not\in M$.
We claim that
\begin{center}
($\ddag$)  $T_jT_k$ divides $T^s$ for some $(j,k)\in\Delta$.
\end{center}

We prove ($\ddag$) by ``induction on $s-q$''.  If $s=q$, then from
(*) we know that either $a_s\not\in I_j$ for some $j$ or ($\ddag$) holds.
But $I_j\subseteq M$ by unshortenability,  so $a_s\not\in I_j$ by the
definition of $s$,  so ($\ddag$) holds.   Now suppose that $\deg(s)>\deg(q)$.
By (*), either $a_q\in I_j$ for some $j$, or $T_jT_k$ divides $T^q$
for some $(j,k)\in\Delta$.   Accordingly there are two cases.

First assume that $a_q\in I_j$ for some $j$.   Then there exists a linear
form $F'=a_qT_j-a_1T_1-\cdots-a_{j-1}T_{j-1}$ in $J$ with $\deg(a_q)
+\deg(x_j)=\deg(a_1)+\deg(x_1)=\cdots=\deg(a_{j-1})+\deg(x_{j-1})$.
By our assumption that $\deg(x_1)\leq\ldots\leq\deg(x_n)$, we have
$\deg(a_1)\geq\ldots\geq\deg(a_{j-1})\geq\deg(a_q)$.    Since $a_q\in M$
by assumption, it follows that $a_1,\ldots,a_{j-1}$ are also in $M$.
Consider $G=F-F'T^q/T_j$.   Since $F'_*=a_qT_j$, it follows that
$F_*=(F'T^q/T_j)_*$, so that $\deg(G_*)<\deg(F_*)$.   Since $F'T^q/T_j$
belongs to $M[T_1,\ldots,T_n]$, it follows that $G_s$ is the term 
of least degree in $G$ with coefficient not in $M$.   By induction, we 
have $T_jT_k$ divides $T^s$ for some $(j,k)$ in $\Delta$.   We are done.

Now assume that $T_jT_k$ divides $T^q$ for some $(j,k)$ in $\Delta$.
By the definition of $\Delta$ and axiom~2 in the definition of quadratic sequence,
there exists a $2$-form $F'$ in $J$ with $F_*'=T_jT_k$.   Now consider
$G=F-a_qF'T^q/T_jT_k$,  and argue just as in the first case.

\bex
{\em (Defining ideals of monomial projective space curves
lying on the quadric $xw-yz=0$ as in [MS])  
This example is treated in \S2 of [RS] (see pages 558, 559)
and in \S3.3 of [RV].   We have
$$
k[\Delta]=k[T_1,\ldots,T_{b-c+2}]/(T_3,\ldots,T_{b-c+1})^2
$$
so that $k[\Delta]$ is Cohen-Macaulay and its Hilbert series is given by
$$
H(k[\Delta];t)=\frac{1+(b-c-1)t}{(1-t)^3}.
$$}
\eex

\bex
{\em (Huckaba-Huneke ideals of analytic deviation $1$ and $2$.)
This example is treated  in 
[RV]   (see \S3.4 of that paper).   We have
$$
k[\Delta]=k[T_1,\ldots,T_n]/(T_{m+1},\ldots,T_{n})^2
$$
so that $k[\Delta]$ is Cohen-Macaulay and its Hilbert series is given by
$$
H(k[\Delta];t)=\frac{1+(n-m)t}{(1-t)^m}
$$
Here $m$ is the analytic spread of the ideal and $n$ is the
minimal number of generators of the ideal.}
\eex

\bex
{\em (Straightening-closed ideals in graded algebras with
straightening law.)
 These examples are treated in \S2 of [RS] and in \S3.2 of [RV].
Let $R$ be a graded algebra with straightening law on a finite poset
$\Pi$ over a field $k$ and let $M$ denote the graded maximal ideal of $R$.
If $\Omega$ is a straightening-closed ideal of $\Pi$ which admits a
linearization $\#:\Omega\to\{1,\ldots,n\}$ satisfying
$$
\deg(x_1)\leq\cdots\leq\deg(x_n),
$$
then $k[\Delta]$ is the face ring of the poset $\Omega$.
Thus $F(I)$ has the same Hilbert series as the face ring of $\Omega$.
And $F(I)$ is Cohen-Macaulay if $\Omega$ is a Cohen-Macaulay poset.
While these conclusions are well-known for most particular examples
of straightening-closed ideals, what appears to be new is
the conclusion from Remark~\ref{remark}
that the reduction number is independent of the minimal reduction
and can be read off the Hilbert series of $F(I)$.   We illustrate
this last conclusion by means of an example:
  Let $(Z_{ij}|1\leq i \leq 2, 1\leq j\leq n)$
be a generic $2\times n$ matrix ($n\geq3$),
$R=k[Z_{ij}]$ the polynomial ring over
a field $k$ in these $2n$ indeterminates, $\Omega$ the poset of $2\times2$
minors of the matrix $Z$,  and $I$ the ideal of $R$ generated by the
elements of $\Omega$.   
The Hilbert series of $F(I)$ is (see \cite{chv})
$$
H(F(I),t) = \frac{h_0+h_1t+\cdots+h_{n-3}t^{n-3}}{(1-t)^{2n-3}}
$$
where
$$
h_i={{n-2}\choose i}^2-{{n-3}\choose{i-1}}{{n-1}\choose{i+1}}.
$$
Thus the reduction number of $I$ is $n-3$.}
\eex

\section{Examples}
In this section we will use the fact that if $I$ is generated by 
homogeneous polynomials of equal degree in a polynomial ring
over a field $k$, then the fiber cone of $I$ is  isomorphic to $k[I]$.

\bex {\em
\label{}
\cite[Example~10.27]{harris} Let $R=k[x,y,z]_m$, where $k$ is a field and
$m=(x,y,z)$.  Let 
$I=((x^{2}, y^{2})^{3}, x(x^{2}, y^{2})z^{3}, z^{6})$. 
Then $F(I)$ is CM, $r(I)=2$ and 
$$
H(F(I), t) =
\f{1+ 4t + t^2}{(1-t)^3} .
$$
We prove by induction on $n$, that 
$$
I^{n}=((x^{2}, y^{2})^{3i} z^{6(n-i)}; i=0 \ldots, n)
      + (x(x^{2}, y^{2})^{3i+1}z^{6(n-1-i)+3}; i=0, \ldots, n-1 ).
$$
For $n=1$, it is easy to verify. Assume that $n \geq 1$. Then
\beqn
&&  I^{n+1}\\
&=& ((x^{2}, y^{2})^{3(i+1)} z^{6(n-i)},
     (x^{2}, y^{2})^{3i} z^{6(n+1-i)},
     x(x^{2}, y^{2})^{3(i+1)} z^{6(n-i)+3}; i=0, \ldots n)\\
&&+  (x(x^{2}, y^{2})^{3i+4}z^{6(n-1-i)+3},
    x(x^{2}, y^{2})^{3i+1}z^{6(n-i)+3},
      x^2(x^{2}, y^{2})^{3i+2}z^{6(n-i)}; i=0, \ldots n-1)\\
&=& ((x^{2}, y^{2})^{3i} z^{6(n+1-i)}; i=0 \ldots, n+1)
  + (x(x^{2}, y^{2})^{3i+2}z^{6(n-i)+3}; i=0, \ldots, n )
\eeqn
since $x^2(x^{2}, y^{2})^{3i+2}z^{6(n-i)} \seq (x^{2}, y^{2})^{3i+3}
z^{6(n-i)}$ for all $i=0, \ldots, n-1$.  
Therefore
\beqn
\mu(I^n) 
&=& \sum_{i=0}^{n} (3i+1) + \sum_{i=0}^{n-1} (3i+2)\\
&=& \sum_{i=0}^{n}(6n+3)  -(3n+2) \\
&=& 6 {n+1 \choose 2} + 1\\
&=& 6 {n+2 \choose 2} - 6 {n+1 \choose 1} 
   + 1.
\eeqn
Hence
\beqn
\sum_{i=0}^{\infty} \mu(I^n) t^n 
= \f{6}{(1-t)^3} - \f{6}{(1-t)^2} + \f{1}{(1-t)}
= \f{1+ 4 t + t^2}{(1-t)^3}.
\eeqn
Let $J=(x^{6} , y^{6}, z^{6})$. We will show that $JI^2 = I^3$. Now, 
\beqn
I^2
&=& ((x^2, y^2)^{6}, x  (x^2, y^2)^{4} z^{3}, 
    (x^2, y^2)^{3} z^{6}, x ( x^2 , y^2) z^{9}, z^{12})\\
I^3 
&=& ((x^{2},y^{2})^{9}, \,  x (x^{2},y^{2})^{7} z^{3}, \, 
    (x^{2},y^{2})^{6} z^{6}, \,  x (x^{2},y^{2})^{4} z^{9}, \,
    (x^{2},y^{2})^{3} z^{12}, \, x (x^{2},y^{2}) z^{15}, \,
     z^{18})\\
&=& z^{6} I^{2} 
  + (x^{6}, y^{6})((x^{2}, y^{2})^{6}, x (x^{2}, y^{2})^{4} z^{3})\\
&\seq& JI^2.
\eeqn
Hence $JI^2 = I^3$,  $J$ is a minimal reduction of
$I$ and $\ell (I/J+mI) = \mu(I) - \mu(J) = 4$. We will show that 
$\ell (I^2/ JI + mI^2 ) = 1$. Now, 
\beqn
JI 
&=& (x^{6} , y^{6}, z^{6})
   ((x^{2}, y^{2})^3, xz^{3}(x^{2}, y^{2}), z^{6}) \\
&=& ((x^2, y^2)^{6}, (x^{7}, x y^{6})  (x^2, y^2) z^{3}, 
    (x^2, y^2)^{3} z^{6}, x( x^2 , y^2) z^{9}, z^{12}) \\
mI^2 
&=& ((x,y)^{13},       (x^2, y^2)^{6} z, 
    x(x,y)^{9} z^3 ,  x(x^2, y^2)^{4} z^{4}, 
   (x,y)^{7} z^{6},  (x^2, y^2)^{3} z^{7},\\
&&  x(x,y)^{3} z^{9}, x( x^2 , y^2) z^{10}, 
     (x,y) z^{12},      z^{13})       \\
JI + mI^2
&=& JI + xm (x^4 y^4) z^3 .
\eeqn    
Hence 
\beqn
    \ell \left( \f{I^2}{JI + mI^2} \right)
&=& \ell \left( \f{R}{JI + mI^2} \right)
  - \ell \left( \f{R}{I^2} \right)\\
&=& \left[ \sum_{i=0}^{11} {n+2 \choose 2} 
  +   73 + 48  + 20 + 4 + 2 \right] \\
&&  -  \left[ \sum_{i=0}^{11} {n+2 \choose 2} 
  +  72 + 48  + 20 + 4 + 2 \right]\\
&=& 1.
\eeqn
In view of Theorem~\ref{main}, $F(I)$ is CM.}
\eex

\bex{\em
\label{hoc}
 \cite[Remark~1]{hoc}
Let $R=k[x_1, x_2, y_1, \ldots, y_s]$ $(s \geq 2)$ and let
$I=(x_iy_j: i=1,2; 1 \leq j \leq s)$ be an ideal of $R$. Then
$F(I) \cong k[x_iy_j: i=1,2; 1 \leq j \leq s] \cong k[T_{ij}: 1
\leq i \leq 2; 1 \leq j \leq s]/I_2$, where $I_2$ is the ideal
generated by the $2 \times 2$ minors of the matrix $(T_{ij})$
\cite[Remark~1]{hoc}.  It is known that $F(I)$ is CM
\cite[Theorem~1$^{0}$]{hoc}..  We show that $F(I)$ is CM by using
Theorem~\ref{main}. First we will show that $r(I) \leq 1$.  
Put 
$$
J= (x_1y_i+x_2y_{i+1}: 1 \leq i \leq s-1; x_1y_s, x_2y_1).
$$
We will prove that $JI = I^2$. Now 
$I^2 = (x_1y_i\,x_1y_j, x_1y_i \, x_2y_j, x_2y_i \, x_2y_j
                                       ; 1 \leq i \leq j \leq s)$. 
Obviously $x_1y_i\, x_1y_s, \in JI$ for $1 \leq i \leq s$. Let $j <
s$. Then   
\beqn
    x_1y_1 \, x_1y_j
&=& x_1y_1 (x_1y_j + x_2y_{j+1}) - x_2y_1 \, x_1y_{j+1} \in JI.
\eeqn
If $i \geq 2$, then by induction hypothesis, 
\beqn
    x_1y_i \, x_1y_j
&=& x_1y_i (x_1y_j + x_2 y_{j+1}) 
  - x_1y_{j+1}(x_1y_{i-1} + x_2y_i)
  + x_1y_{j+1} \, x_1y_{i-1} \in JI.
\eeqn
Now,
$$
x_1 y_1 \, x_2 y_j = x_1 y_j \, x_2 y_1 \in JI
$$
for  $1 \leq j  \leq s$. If $i \geq 2$, then
\beqn
x_1y_i \, x_2y_j 
&=& x_1y_i (x_1y_{j-1} + x_2y_j) - x_1y_{i} \, x_1y_{j-1} \in JI
\eeqn
for $1 \leq j \leq s$. 
Since $x_2 y_1 \in J$, $x_2 y_1 \, x_2 y_j \in JI$ for $1 \leq j \leq
s$. Hence for $2 \leq 
    i \leq j \leq s$,
\beqn
    x_2y_i \, x_2y_j
&=& x_2y_i (x_1y_{j-1} + x_2y_j) - x_2y_i \, x_1y_{j-1} \in JI.
\eeqn
This shows that $JI = I^2$. We will now compute $H(F(I),t)$. 
Note that  $I^n$ is generated by monomials $x_1^{h_1} x_2^{h_2}
y_1^{l_1} \ldots y_s^{l_s} $ where $h_1+h_2=n$ and
$l_1 + \cdots + l_s = n$. Hence,
\beqn
\sum_{n \geq 0} \mu(I^n) t^n
&=& \sum_{n \geq 0} \left[ {n+1 \choose 1}{n+s-1 \choose s-1}
                                               \right] t^n\\
&=& \sum_{n \geq 0} \left[ s {n+s-1 \choose s}
    + {n+s-1 \choose s-1 } \right] t^n\\
&=& \sum_{n \geq 0} \left[ s {n+s \choose s} 
                     - (s-1) {n+s-1 \choose s-1} \right] t^n\\
&=& \f{s - (s-1)(1-t)}{(1-t)^{s+1}}\\
&=& \f{1 + (s-1)t}{(1-t)^{s+1}}.
\eeqn
This shows that $a(I) = s+ 1$ and hence $J$ is a minimal reduction of
$I$. 
Since $e(F(I)) = s =
\mu(I) - \dim \, F(I) + 1$, $F(I)$ is CM with  minimal multiplicity. }
\eex

\end{document}